\documentclass[12pt]{amsart}
\usepackage{amssymb}
\usepackage[all]{xy}

\newcommand{\issuenumber}{12}
\newcommand{\issuemonth}{March}
\newcommand{\issueyear}{2005}

\newtheorem{thm}{Theorem}[section]
\newtheorem{prob}[thm]{Problem}

\newtheorem{issue}{Issue}

\theoremstyle{definition}

\theoremstyle{remark}

\newcommand{\alphs}{{\aleph_0}}
\newcommand{\ed}{\end{thebibliography}\general\end{document}}

\renewcommand{\>}{\right >}

\renewcommand{\b}{\mathfrak{b}}

\newcommand{\p}{\mathfrak{p}}

\newcommand{\NON}{{\mathsf   {NON}}}
\newcommand{\COF}{{\mathsf   {COF}}}

\newcommand{\dnannouncement}[1]{[\S\ref{#1} below]}

\newcommand{\M}{\mathcal{M}}

\newcommand{\cov}{\mathsf{cov}}

\newcommand{\cf}{\mathsf{cf}}

\newcommand{\R}{\mathbb{R}}
\newcommand{\Q}{\mathbb{Q}}

\newcommand{\fo}{\mathfrak{od}}

\renewcommand{\b}{\mathfrak{b}}

\renewcommand{\split}{\mathsf{Split}}
\newcommand{\bq}{\begin{quote}}
\newcommand{\eq}{\end{quote}}

\newcommand{\B}{\mathcal{B}}
\newcommand{\BG}{\B_\Gamma}

\newcommand{\sone}{\mathsf{S}_1}    \newcommand{\sfin}{\mathsf{S}_{fin}}
\newcommand{\ufin}{\mathsf{U}_{fin}}

\newcommand{\nin}{\not\in}

\newcommand{\naturals}{{\mathbb N}}
\newcommand{\N}{\naturals}

\newcommand{\sbst}{\subseteq}
\newcommand{\by}[2]{\par\hfill\emph{#1}, #2}
\newcommand{\nby}[1]{\par\hfill\emph{#1}}
\newcommand{\Tau}{\mathrm{T}}
\newcommand{\CE}{\textsc{CE}}

\newcommand{\be}{\begin{enumerate}}
\newcommand{\ee}{\end{enumerate}}
\newcommand{\bi}{\begin{itemize}}
\newcommand{\ei}{\end{itemize}}
\renewcommand{\i}{\item}

\newcommand{\general}{\small\vfill\par\noindent\hrulefill\par
\noindent\textbf{Previous issues.} The first issues of this
bulletin, which contain general information (first issue), basic
definitions, research announcements, and open problems (all
issues) are available online, on \arx{math.GN/$x$}, where $x$ is
\texttt{0301011}, \texttt{0302062}, \texttt{0303057},
\texttt{0304087}, \texttt{0305367}, \texttt{0312140},
\texttt{0401155}, \texttt{0403369}, \texttt{0406411},
\texttt{0409072}, and \texttt{0412305},
respectively, for issues number $1$ to $11$.\\[0.1cm]
\textbf{Contributions.}
Please submit your contributions (announcements, discussions, and open problems)
by e-mailing us. It is preferred to write them
in \LaTeX{}.
The authors are urged to use as standard notation as possible, or otherwise give
the definitions or a reference to where the notation is explained.
Contributions to this bulletin would not require any transfer of copyright,
and material presented here can be published elsewhere.\\[0.1cm]
\textbf{Subscription.}
To receive this bulletin (free) to your
e-mailbox, e-mail us:\\
{boaz.tsaban@weizmann.ac.il}
}

\newcommand{\nArxPaper}[5]{\subsection{#2}{#4}\par\hfill{\arx{#1}}\par\hfill\emph{#3}}

\newcommand{\nAMSPaper}[5]{\subsection{#2}{#4}\par\hfill{\texttt{#1}}\par\hfill\emph{#3}}
\newcommand{\SPMBul}{\textbf{$\mathcal{SPM}$ Bulletin}}

\newcommand{\arx}[1]{\texttt{http://arxiv.org/abs/#1}}
\newcommand{\url}[1]{\bq\texttt{#1}\eq}
\newcommand{\online}[1]{The paper is available online at \url{#1}}

\title[$\mathcal{SPM}$ Bulletin \textbf{\issuenumber} (\issuemonth{} \issueyear)]{%
$\mathcal{SPM}$ Bulletin\\[0.5cm]
Issue number \issuenumber: \issuemonth{} \issueyear{} \CE{}}

\begin{document}
\maketitle

\tableofcontents

\section{Editor's note}

\textbf{1. Proceedings of first SPM Workshop.}
The proceedings of the first workshop on \emph{Coverings, Selections, and Games in Topology}
(Lecce, Italy, June 27--29, 2002) were published as a special issue of \textbf{Note di Matematica}
(Volume 22, Issue 2 (2003)), and are now also available online at:
\bq
\verb|http://siba2.unile.it/ejournals/search/issue.php?single_pub=1&|\\
\verb|pubid=1&id=342&barcol=000033&recs=10|
\eq
(Alternatively, go to \texttt{http://siba2.unile.it/notemat/},
choose ``Contents, Abstracts, Full text articles'', then under the year 2003, choose
``vol.\ 22, issue 2'').
In \dnannouncement{SPMproc} we quote the editors' preface to the volume.

The next SPM Workshop is planned for the coming December,
details are yet to come.

\medskip\textbf{2. Mathematical breakthroughs.} Justin Moore announces in the Mathematics ArXiv a sequence of
astonishing results. Liljana Babinkostova has solved an 1938 problem of Rothberger, and established a game
theoretic characterization of countable dimensionality.
These and other beautiful mathematical results are announced below.

\medskip\textbf{3. Mathematics ArXiv.}
Elliot Peal, editor of the \emph{Topology Atlas}, recommends storing preprints in the ArXiv.
We quote his message from the bulletin \emph{Topology News}:
\bq
Topology Atlas maintains a preprint server. Our preprint server was most
active in 1996 and 1997 and nearly all of our preprints have been
published elsewhere by now. We suggest that you use the mathematics arXiv
(arXiv.org) or its Front (front.math.ucdavis.edu) for storing and finding
preprints. The arXiv is a powerful and stable resource for distributing
results in active research communities. Upon request, we can help you
archive all your TeX files in the arXiv.
\eq
We support this recommendation, and also recommend that you get subscribed
to the mailing lists of the arXiv, the categories of interest to readers
of the \SPMBul{} are probably those of General Topology and Logic.
Details on subscription etc.\ are available at
\url{http://arxiv.org/}
\medskip

\medskip\textbf{3. Change of name.} We were informed by an anonymous referee
that the name $\mathfrak{o}$ used in Issue 10 of the \SPMBul{} (and its subsequent
issues) for the \emph{$o$-diagonalization} number is already used for
Leathrum's \emph{off-branching} number \cite{leathrum}.
We therefore move to use $\mathfrak{od}$
for the new cardinal.

\medskip

Contributions to the next issue are, as always, welcome.

\medskip

\by{Boaz Tsaban}{boaz.tsaban@weizmann.ac.il}

\hfill \texttt{http://www.cs.biu.ac.il/\~{}tsaban}

\section{Proceedings of the first workshop on Coverings, Selections, and Games in Topology}\label{SPMproc}

The first Workshop on Coverings, Selections and Games in Topology
was held in Lecce, Italy, on June 27{29, 2002 and organized by the Department
of Mathematics of Lecce University.

The aim of this workshop was to consider development and activities over
the last few years in a quickly growing field of Mathematics known under the
name Selection Principles and to discuss its relationships with other areas of
Mathematics, such as (generalized) Ramsey theory, infinite game theory, combinatorial
cardinals, function spaces, hyperspaces, uniform spaces, topological
algebras and so on.

The work was organized in such a way as to have plenary 50 minute lectures
on general aspects as well as on some specific topics and 25 minute contributed
talks. Each session was ended by a discussion on presented results.
The conference has been closed by a Round Table. The goal was to consider
and discuss further common research activities, open problems in the field and
the organization of new future meetings on the same topic. It was decided to
publish the Proceedings of the meeting as a special issue of Note di Matematica
following the usual referee procedure established by the journal. In particular,
it has been decided the Proceedings to contain two survey papers: the first
one opens this volume and presents general aspects of the theory of selection
principles with recent progress in this area and with directions of further investigation;
the second one, which closes the volume, is devoted to the presentation
and discussion of open problems. Special thanks are due to the authors of these
surveys.

Also, it was planned to edit the SPM Bulletin, an electronic bulletin dedicated
to the field with announcements of new results as well as with open
problems.

We hope that these Proceedings which contain some of the papers presented
at the Workshop will help to attract new researchers and stimulate the investigation
in this field of Mathematics.
\nby{Cosimo Guido, Ljubisa D.\ R.\ Ko\v{c}inac, and Marion Scheepers}

\subsection{Contents of the proceedings volume}
~\be
\i C.\ Guido,  L.\ D.\ R.\ Ko\v{c}inac, M.\ Scheepers:\\
\emph{Preface} (p.~1).
\i M.\ Scheepers:\\
\emph{Selection principles and covering properties in Topology} (p.~3).
\i M.\ Sakai:\\
\emph{The Pytkeev property and the Reznichenko property in function spaces} (p.~43).
\i B.\ Tsaban:\\
\emph{Selection principles and the Minimal Tower problem} (p.~53).
\i T.\ Weiss, B.\ Tsaban:\\
\emph{Topological diagonalizations and Hausdorff dimension} (p.~83).
\i D.\ Leseberg:\\
\emph{Symmetrical extensions and generalized nearness} (p.~93).
\i M.\ Caldas, D.\ N.\ Georgiou, S.\ Jafari, T.\ Noiri:\\
\emph{More on $\delta$-semiopen sets } (p.~113).
\i L.\ D.\ R.\ Kocinac:\\
\emph{Selection principles in uniform spaces} (p.~127).
\i S.\ D.\ Iliadis:\\
\emph{Saturated classes of bases} (p.~141).
\i S.\ E.\ Han:\\
\emph{Generalized digital $(k_0,k_1)$-homeomorphism} (p.~157).
\i L.\ Babinkostova, M.\ Scheepers:\\
\emph{Combinatorics of open covers (IX): Basis properties} (p.~167).
\i B.\ Tsaban:\\
\emph{Selection principles in mathematics: A milestone of open problems} (p.~179).
\ee

\section{Additional research announcements}

\nArxPaper{math.LO/0501525}
{A five element basis for the uncountable linear orders}
{Justin Tatch Moore}
{In this paper I will show that it is relatively consistent with the usual
axioms of mathematics (ZFC) together with a strong form of the axiom of
infinity (the existence of a supercompact cardinal) that the class of
uncountable linear orders has a five element basis. In fact such a basis
follows from the Proper Forcing Axiom, a strong form of the Baire Category
Theorem. The elements are $X$, $\omega_1$, $\omega_1^*$, $C$, $C^*$ where $X$ is any suborder
of the reals of cardinality $\aleph_1$ and $C$ is any Countryman line. This confirms
a longstanding conjecture of Shelah.}

\nArxPaper{math.LO/0501526}
{Set mapping reflection}
{Justin Tatch Moore}
{In this note we will discuss a new reflection principle which follows from
the Proper Forcing Axiom. The immediate purpose will be to prove that the
bounded form of the Proper Forcing Axiom implies both that $2^\alphs = \aleph_2$
and that $L(P(\aleph_1))$ satisfies the Axiom of Choice. It will also be
demonstrated that this reflection principle implies that combinatorial
principle $\Box(\kappa)$ fails for all regular $\kappa > \aleph_1$.
}

\nArxPaper{math.LO/0501527}
{The Proper Forcing Axiom, Prikry forcing, and the Singular Cardinals Hypothesis}
{Justin Tatch Moore}
{The purpose of this paper is to present some results which suggest that the
Singular Cardinals Hypothesis follows from the Proper Forcing Axiom. What will
be proved is that a form of simultaneous reflection follows from the Set
Mapping Reflection Principle, a consequence of PFA. While the results fall
short of showing that MRP implies SCH, it will be shown that MRP implies that
if SCH fails first at kappa then every stationary subset of $S_{\kappa^+}^\omega =
\{\alpha < \kappa^+ : \cf(\alpha) = \omega\}$ reflects. It will also be demonstrated that MRP
always fails in a generic extension by Prikry forcing.
}

\nArxPaper{math.GN/0501524}
{A solution to the $L$ space problem and related ZFC constructions}
{Justin Tatch Moore}
{In this paper I will construct a non-separable hereditarily Lindelof space ($L$
space) without any additional axiomatic assumptions. I will also show that
there is a function $f:[\omega_1]^2\to \omega_1$ such that if $A,B\sbst\omega_1$,
are uncountable and $\xi < \omega_1$, then there are $\alpha < \beta$ in $A$
and $B$ respectively with $f\{\alpha,\beta\} = \xi$.

Previously it was unknown whether such a function existed even if $\omega_1$ was replaced
by $2$.
Finally, I will prove that there is
no basis for the uncountable regular Hausdorff spaces of cardinality $\aleph_1$.
Each of these results gives a strong refutation of a well known and
longstanding conjecture. The results all stem from the analysis of oscillations
of coherent sequences $\<e_\alpha : \alpha < \omega_1\>$  of finite-to-one functions. I expect
that the methods presented will have other applications as well.
}

\nArxPaper{math.GN/0501076}
{Countable Tightness, Elementary Submodels and Homogeneity}
{Ramiro de la Vega}
{We show (in ZFC) that the cardinality of a compact homogeneous space of
countable tightness is no more than the size of the continuum.}

\subsection{No transcendence basis of $\R$ over $\Q$ can be an analytic set}
It has been proved by Sierpi\'nski that no linear basis of $\R$ over $\Q$
can be an analytic set. Here we show that the same assertion holds
by replacing ``linear basis'' with ``transcendence basis''. Furthermore,
it is demonstrated that purely transcendental subfields of $\R$ generated
by Borel bases of the same cardinality are Borel isomorphic (as fields).
Following Mauldin's arguments, we also indicate, for each ordinal $\alpha$ such
that $1\le\alpha < \aleph_1$ ($2\le\alpha<\aleph_1$), the existence of subfields
of $\R$ of exactly additive (multiplicative, ambiguous) class $\alpha$ in $\R$.
\nby{Enrico Zoli}

\subsection{Two properties of $C_p(X)$ weaker than Fr\'echet Urysohn property}
For a Tychonoff space $X$,
we denote by $C_p(X)$
the space of all real-valued continuous functions on $X$ with
the topology of pointwise convergence.
In this paper,
we study $\kappa$-Fr\'echet Urysohn property
and weak Fr\'echet Urysohn property of $C_p(X)$.
Our main results are that:
\be
\i $C_p(X)$ is $\kappa$-Fr\'echet Urysohn iff
$X$ has property ($\kappa$-FU) (i.e.
every pairwise disjoint sequence of finite subsets of $X$
has a strongly point-finite subsequence),
in particular a Baire space $C_p(X)$ is $\kappa$-Fr\'echet Urysohn;
\i among separable metrizable spaces, every $\lambda$-space has property
($\kappa$-FU) and
every space having property ($\kappa$-FU) is always of the first category;
\i every analytic space has the $\omega$-grouping property,
hence
for every analytic space $X$,
$C_p(X)$ is weakly Fr\'echet Urysohn.
\ee
\nby{Masami Sakai}

\nArxPaper{math.LO/0501421}
{Some partition properties for measurable colourings of $(\aleph_1)^2$}
{James Hirschorn}
{We construct a measure on $(\aleph_1)^2$ over the ground model in the forcing
extension of a measure algebra, and investigate when measure theoretic
properties of some measurable colouring of $(\aleph_1)^2$ imply the existence of
an uncountable subset of $\aleph_1$ whose square is homogeneous. This gives a
new proof of the fact that, under a suitable axiomatic assumption, there are no
Souslin $(\aleph_1,\aleph_1)$ gaps in the Boolean algebra $L^0(\nu)/Fin$ when $\nu$
is a separable measure.

To appear in: \emph{Proceedings of the Kyoto conference on Forcing Method and Large
Cardinals, 2004}.
}

\nArxPaper{math.LO/0502394}
{Potential theory and forcing}
{Jindrich Zapletal}
{We isolate a combinatorial property of capacities leading to a construction
of proper forcings. Then we show that many classical capacities such as the
Newtonian capacity satisfy the property.}

\nAMSPaper{http://www.ams.org/journal-getitem?pii=S0002-9939-05-07799-3}
{On decompositions of Banach spaces
of continuous functions on Mr\'{o}\-wka's spaces}
{Piotr Koszmider}
{It is well known that if
$K$ is infinite compact Hausdorff and scattered
(i.e., with no perfect subsets),
then the Banach space $C(K)$
of continuous functions on $K$ has
complemented copies of $c_{0}$,
i.e., $C(K)\sim c_{0}
\oplus X\sim c_{0}\oplus c_{0} \oplus X\sim c_{0}\oplus C(K)$.
We address the question if this could be the
only type of decompositions of $C(K)\not \sim c_{0}$
into infinite-dimensional summands for $K$ infinite, scattered.
Making a special set-theoretic assumption such as the
continuum hypothesis or Martin's axiom we construct an example
of Mr\'{o}wka's space (i.e., obtained from an
almost disjoint family of sets of positive integers) which answers
positively the above question.}

\nArxPaper{math.GN/0503275}
{A note on $D$-spaces}
{Gary Gruenhage}
{We introduce notions of nearly good relations and
$N$-sticky modulo a relation as tools for proving  that spaces are
$D$-spaces. As a corollary to general results about such
relations, we show that $C_p(X)$ is hereditarily a $D$-space
whenever $X$ is a Lindel\"of $\Sigma$-space. This answers a
question of Matveev, and improves a result of Buzyakova, who
proved the same result for $X$ compact.

We also prove that if a space $X$ is the union of finitely many
$D$-spaces, and has countable extent, then $X$ is linearly Lindel\"of.
It follows that if $X$ is in addition countably compact, then $X$ must
be compact.  We also show that Corson compact spaces are
hereditarily $D$-spaces.  These last two results answer recent
questions of Arhangel'skii.  Finally, we answer a question of van
Douwen by showing that a perfectly normal collectionwise-normal
non-paracompact space constructed by R.\ Pol is a $D$-space.}

\subsection{Set-theoretic properties of Schmidt's ideal}
We study some set-theoretic properties of Schmidt's $\sigma$-ideal on $\R$,
emphasizing its analogies and dissimilarities with both the classical $\sigma$-ideals on
$\R$ of Lebesgue measure zero sets and of Baire first category sets. We highlight
the strict analogy between Schmidt's ideal on $\R$ and Mycielski's ideal on $\{0,1\}^\N$.
\nby{Marcin Kysiak and Enrico Zoli}

\nAMSPaper{http://www.ams.org/journal-getitem?pii=S0002-9939-05-07824-X}
{Almost-disjoint coding and strongly saturated ideals}
{Paul B.\ Larson}
{We show that Martin's Axiom plus $\mathfrak{c} = \aleph_{2}$
implies that there is no $(\aleph_{2},
\aleph_{2},\aleph_{0})$-saturated $\sigma$-ideal on $\omega_{1}$.
}

\section{Selective screenability and covering dimension}
Let $X$ be a topological space. In \cite{bing} Bing introduced the following notion of \emph{screenability}: For each open cover $\mathcal{U}$ of $X$ there is a sequence $(\mathcal{V}_n:n<\infty)$ such that: For each $n$, $\mathcal{V}_n$ is a family of pairwise disjoint open sets; for each $n$, $\mathcal{V}_n$ refines $\mathcal{U}$ and $\cup_{n<\infty}\mathcal{V}_n$ is an open cover of $X$. In \cite{AG} Addis and Gresham introduced the selective version screenability, defined as follows: For each sequence $(\mathcal{U}_n:n<\infty)$ of open covers of $X$ there is a sequence $(\mathcal{V}_n:n<\infty)$ such that: For each $n$, $\mathcal{V}_n$ is a family of pairwise disjoint open sets; for each $n$, $\mathcal{V}_n$ refines $\mathcal{U}_n$ and $\cup_{n<\infty}\mathcal{V}_n$ is an open cover of $X$. It is evident that selective screenability implies screenability.

Selective screenability is an example of the following selection principle which was introduced in \cite{babinkostova}: Let $S$ be a set and let $\mathcal{A}$ and $\mathcal{B}$ be families of collections of subsets of the set $S$. Then ${\sf S}_c(\mathcal{A},\mathcal{B})$ denotes the statement that for each sequence $(\mathcal{U}_n:n<\infty)$ of elements of $\mathcal{A}$ there is a sequence $(\mathcal{V}_n:n<\infty)$ such that
\begin{enumerate}
\item{For each $n$, $\mathcal{V}_n$ is a family of pairwise disjoint sets;}
\item{For each $n$, $\mathcal{V}_n$ refines $\mathcal{U}_n$ and}
\item{$\cup_{n<\infty}\mathcal{V}_n$ is a member of $\mathcal{B}$.}
\end{enumerate}
\newcommand{\sC}{{\sf S}_c}
\newcommand{\gC}{{\sf G}_c}
With $\mathcal{O}$ denoting the collection of all open covers of topological space $X$, $\sC(\mathcal{O},\mathcal{O})$ is selective screenability.

Addis and Gresham noted that countable dimensional metrizable spaces are selectively screenable, and asked if the converse is true. Pol showed in \cite{Pol} that the answer is no. The author showed that the countable dimensional metric spaces are exactly characterized by a game-theoretic version of selective screenability.
The following game, denoted $\gC(\mathcal{A},\mathcal{B})$, is naturally associated with $\sC(\mathcal{A},\mathcal{B})$: Players ONE and TWO play as follows: In the $n$-th inning ONE first chooses $\mathcal{O}_n$, a member of $\mathcal{A}$, and then TWO responds with $\mathcal{T}_n$ which is pairwise disjoint and  refines $\mathcal{O}_n$. A play $(\mathcal{O}_1, \mathcal{T}_1, \cdots, \mathcal{O}_n, \mathcal{T}_n, \cdots)$ is won by TWO if $\cup_{n<\infty}\mathcal{T}_n$ is a member of $\mathcal{B}$; else, ONE wins. We can consider versions of different length of this game as follows: For an ordinal number $k$ let $\gC^k(\mathcal{A},\mathcal{B})$ be the game played as follows: in the $l$-th inning ($l<k$) ONE first chooses $\mathcal{O}_l$, a member of $\mathcal{A}$, and then TWO responds with a pairwise disjoint $\mathcal{T}_l$ which refines $\mathcal{O}_l$. A play
\[
  \mathcal{O}_0,\mathcal{T}_0, \cdots, \mathcal{O}_{\l}, \mathcal{T}_{\l},\cdots \, l<k
\]
is won by TWO if $\cup_{l<k}\mathcal{T}_l$ is a member of $\mathcal{B}$; else, ONE wins. Thus the game $\gC(\mathcal{A},\mathcal{B})$ is $\gC^{\omega}(\mathcal{A},\mathcal{B})$.\\
The author showed the following:
\begin{thm}\label{twoswinningstrategy} Let $X$ be a metric space.
\begin{enumerate}
\item{If $X$ is countable dimensional, then TWO has a winning strategy in $\gC^{\omega}(\mathcal{O},\mathcal{O})$.}
\item{If TWO has a winning strategy in $\gC^{\omega}(\mathcal{O},\mathcal{O})$, then $X$ is countable dimensional.}
\end{enumerate}
\end{thm}

In Pol's example, TWO has a winning strategy in the game $\gC^{\omega+1}(\mathcal{O},\mathcal{O})$.

\begin{thm}\label{findim}Let $X$ be a metric space. The following are equivalent:
\begin{enumerate}
\item{If $X$ is $n$-dimensional then TWO has a winning strategy in $\gC^{n+1}(\mathcal{O},\mathcal{O})$.}
\item{If TWO has a winning strategy in $\gC^{n+1}(\mathcal{O},\mathcal{O})$, then $X$ is $n$-dimensional.}
\end{enumerate}
\end{thm}

From this Theorem we obtain that metric space $X$ is $n$-dimensional if, and only if, TWO has a winning strategy in $\gC^{n+1}(\mathcal{O},\mathcal{O})$ but not in $\gC^{n}(\mathcal{O},\mathcal{O})$.

\nby{Liljana Babinkostova}

\section{On a problem of Rothberger and Sierpinski}
Let $Y$ be a subspace of the metric space $X$. Then $\mathcal{O}_X$ denotes the collection of open covers of $X$, and $\mathcal{O}_{XY}$ denotes the collection of open covers of $Y$ by sets open in $X$. Let $\mathcal{F}_X$ denote the collection of finite open covers of $X$.

In 1924 Menger defined in \cite{MENGER} the following basis property, denoted M by Sierpi\'nski: A metric space $X$ has property M if there is for each basis $\mathcal{B}$ of the space a sequence $(B_n:n<\infty)$ from the base such that the diameters of the $B_n$'s converge to zero, and the $B_n$'s cover $X$.

Hurewicz showed in \cite{HURE25} that a metric space $X$ has Menger's basis property M if, and only if, it has ${\sf S}_{fin}(\mathcal{O}_X,\mathcal{O}_{X})$.

According to Rothberger Sierpi\'nski also defined the basis property M' thus: A metric space $X$ has property M' if there is for each basis $\mathcal{B}$ and each sequence $(\epsilon_n:n<\infty)$ of positive real numbers, a sequence $(B_n:n<\infty)$ of elements of $\mathcal{B}$ such that for each $n$ $diam(B_n)<\epsilon_n$, and $\{B_n:n<\infty\}$ covers $X$.

Let $M'(X,Y)$ denote the version of $M'$ where we require that the sequence of $B_n$'s cover the subspace $Y$. (This is the  relative version of the Rothberger basis property.)

Rothberger showed in Theorem 7 of \cite{Roth38} that ${\sf S}_1(\mathcal{O}_X,\mathcal{O}_X)$ implies $M'$. He posed these problems in
\cite{Roth38}:

\flushleft{\textbf{Problem A:}} Does $M'$ imply ${\sf S}_1(\mathcal{O}_X,\mathcal{O}_X)$?

\flushleft{\textbf{Problem B:}} What is the relationship between $M'$ and ${\sf S}_1(\mathcal{F}_X,\mathcal{O}_X)$?

Fremlin and Miller proved in Theorem 6 of \cite{FM} that $M\not\Rightarrow{\sf S}_1(\mathcal{F}_X,\mathcal{O}_X)$, and that ${\sf S}_1(\mathcal{F}_X,\mathcal{O}_X)\not\Rightarrow M$.

We solve problem A by proving:

{\flushleft \textbf{Theorem A.}} For $X$ a metrizable space with ${\sf S}_{fin}(\mathcal{O}_X,\mathcal{O}_X)$, the following are equivalent:
\begin{enumerate}
\item{${\sf S}_1(\mathcal{O}_X,\mathcal{O}_{XY})$.}
\item{$M'(X,Y)$ holds.}
\end{enumerate}

Note that if $Y=X$, then the hypothesis ${\sf S}_{fin}(\mathcal{O}_X,\mathcal{O}_X)$ is not needed, since by Hurewicz's theorem, it follows from the corresponding properties $M'$ or ${\sf S}_1(\mathcal{O}_X,\mathcal{O}_X)$. And when $Y=X$, Theorem A solves Rothberger's problem, Problem A.

This also gives a solution to Problem B:\\
{\flushleft \textbf{Theorem B.}} $M'\Rightarrow {\sf S}_1(\mathcal{F}_X,\mathcal{O}_{X})$, but
${\sf S}_1(\mathcal{F}_X,\mathcal{O}_{X})\not\Rightarrow M'$.

\nby{Liljana Babinkostova}

\section{Problem of the Issue}

\begin{prob}\label{pvbl}
Is the following statement provable in ZFC?
\bq
Let $M$ be a Baire metric space of weight $\aleph_1$, $A\subset{\mathbb R}$
a perfectly meager set of cardinality $\aleph_1$ and $f:M\to A$ a continuous mapping. Then there is a nonempty open
set $U\subset M$ such that $f$ is constant on $U$.
\eq
\end{prob}

The statement in Problem \ref{pvbl} is true under MA+$\lnot$CH \cite{K}.
A positive answer to Problem \ref{pvbl} yields a ZFC example of a Banach space which is weak
Asplund but whose dual is not weak* fragmentable \cite{K}.

\begin{prob}\label{con}
Is it consistent with ZFC (using large cardinals) that there is a Baire metric space $M$ of weight $\aleph_1$
and a partition $\mathcal U$ of $M$ into meager sets such that:
\be
\i $|\mathcal{U}|=\aleph_1$, and
\i $\bigcup {\mathcal U}'$ has the Baire property in $M$ for each ${\mathcal U}'\subset\mathcal U$?
\ee
\end{prob}

A positive answer to Problem \ref{con} implies a negative anwer to Problem \ref{pvbl}.
Indeed, let $M$ and $\mathcal U$ have the mentioned
properties and let $A$ be any perfectly meager subset of $\mathbb R$ of cardinality $\aleph_1$.
Choose a bijection $\varphi:A\to{\mathcal U}$ and define $g:M\to A$ by $g(m)=a$ if $m\in\varphi(a)$.
Then $g$ has the Baire property and hence there is a residual set $M'\subset M$ with $g_{|M'}$ continuous.
Then $M'$, $A$ and $f=g_{|M'}$ show that the answer to Problem \ref{pvbl} is negative.

A positive answer to Problem \ref{con} implies that $\aleph_1$ is measurable in some transitive model
of ZFC containing all ordinals \cite{FK}.

If there is a precipitous ideal on $\aleph_1$, there is, by \cite{FK}
a Baire metric space of
weight $2^{\aleph_1}$ and a partition of it into meager sets satisfying (1) and (2) of Problem \ref{con}.

\nby{Ondrej Kalenda}

\section{Problems from earlier issues}
In this section we list the still open problems among
the past problems posed in the \SPMBul{}
(in the section \emph{Problem of the month/issue}).
For definitions, motivation and related results, consult the
corresponding issue.

For conciseness, we make the convention that
all spaces in question are
zero-dimentional, separable metrizble spaces.

\begin{issue}
Is $\binom{\Omega}{\Gamma}=\binom{\Omega}{\Tau}$?
\end{issue}

\begin{issue}
Is $\ufin(\Gamma,\Omega)=\sfin(\Gamma,\Omega)$?
And if not, does $\ufin(\Gamma,\Gamma)$ imply
$\sfin(\Gamma,\Omega)$?
\end{issue}

\stepcounter{issue}

\begin{issue}
Does $\sone(\Omega,\Tau)$ imply $\ufin(\Gamma,\Gamma)$?
\end{issue}

\begin{issue}
Is $\p=\p^*$? (See the definition of $\p^*$ in that issue.)
\end{issue}

\begin{issue}
Does there exist (in ZFC) an uncountable set satisfying $\sone(\BG,\B)$?
\end{issue}

\stepcounter{issue}

\begin{issue}
Does $X \nin \NON(\M)$ and $Y\nin\mathsf{D}$ imply that
$X\cup Y\nin \COF(\M)$?
\end{issue}

\begin{issue}
Is $\split(\Lambda,\Lambda)$ preserved under taking finite unions?
\end{issue}
\begin{proof}[Partial solution]
Consistently yes (Zdomsky). Is it ``No'' under CH?
\end{proof}

\begin{issue}
Is $\cov(\M)=\fo$? (See the definition of $\fo$ in that issue.)
\end{issue}

\begin{issue}
Does $\sone(\Gamma,\Gamma)$ always contain an element of cardinality $\b$?
\end{issue}

\begin{thebibliography}{00}

\bibitem{AG}
D.\ F.\ Addis and J.\ Gresham,
\emph{A class of infinite dimensional spaces Part I: Dimension theory and Alexandroff's problem},
Fundamenta Mathematicae \textbf{101} (1978),
195--205.

\bibitem{babinkostova}
L.\ Babinkostova,
\emph{Selection principles in topology} (in Macedonian),
Ph.D.\ thesis (2001).

\bibitem{bing}
R.H.\ Bing,
\emph{Metrization of topological spaces},
Canadian Journal of Mathematics \textbf{3} (1951),
175--186.

\bibitem{FK}
R.\ Frankiewicz and K.\ Kunen,
\emph{Solution of Kuratowski's problem on function having the Baire property, I},
Fundamenta Mathematicae \textbf{128} (1987),
171--180.

\bibitem{FM}
D.\ H.\ Fremlin and A.\ W.\ Miller,
\emph{On some properties of Hurewicz, Menger and Rothberger},
Fundamenta Mathematica \textbf{129} (1988),
17--33.

\bibitem{HURE25}
W.\ Hurewicz,
\emph{\"Uber eine Verallgemeinerung des Borelschen Theorems},
Mathematische Zeitschrift \textbf{24} (1925),
401--421.

\bibitem{K}
O.\ Kalenda,
\emph{A weak Asplund space whose dual is not in Stegall's class},
Proceedings of the American Mathematical Society \textbf{130} (2002),
2139--2143.

\bibitem{leathrum}
T.\ E.\ Leathrum,
\emph{A Special Class of Almost Disjoint Families},
The Journal of Symbolic Logic \textbf{60} (1995),
879--891.

\bibitem{MENGER}
K.\ Menger,
\emph{Einige \"Uberdeckungss\"atze der Punktmengenlehre},
Sitzungsberichte der Wiener Akademie \textbf{133} (1924),
421--444.

\bibitem{Pol}
R.\ Pol,
\emph{A weakly infinite-dimensional compactum which is not countable dimensional},
Proceedings of the American Mathematical Society \textbf{82} (1981),
634--636.

\bibitem{Roth38}
F.\ Rothberger,
\emph{Eine Versch\"arfung der Eigenschaft {\sf C}},
Fundamenta Mathematicae \textbf{30} (1938),
50--55.

\ed